\AtBeginDocument{%
  \paperwidth=\dimexpr
    1in + \oddsidemargin
    + \textwidth
    + 1in + \oddsidemargin
  \relax
  \paperheight=\dimexpr
    1in + \topmargin
    + \headheight + \headsep
    + \textheight
    + 1in + \topmargin
  \relax
  \usepackage[pass]{geometry}\relax
}
\documentclass[smallextended]{svjour3}       
\usepackage{float}
\usepackage[Algorithm,ruled]{algorithm}
\usepackage{algorithmicx}
\usepackage[noend]{algpseudocode}
\usepackage{tabularx}
\newcolumntype{L}{>{\raggedright\arraybackslash}X}
\usepackage{graphicx}
\usepackage{amsmath}
\usepackage{amssymb}
\usepackage{latexsym}
\usepackage{epsfig}
\usepackage[mathscr]{eucal}
\usepackage{graphics} \usepackage{color} \usepackage{multicol}
\usepackage{enumerate} \usepackage{array}

\begin{document}
\title{IP Solutions for International Kidney Exchange Programmes}
\author{P\'eter Bir\'o \thanks{Supported by the Hungarian Academy of Sciences under its Momentum Programme (LP2016-3/2018) and Cooperation of Excellences Grant (KEP-6/2018), and by the Hungarian Scientific Research Fund -- OTKA (no.\ K129086).} \and
M\'arton Gyetvai \and Radu-Stefan Mincu \and
Alexandru Popa\and
Utkarsh Verma}
\institute{P\'eter Bir\'o \and M\'arton Gyetvai \at Institute of Economics, Hungarian Academy of Sciences, Hungary \\
Department of Operations Research and Actuarial Sciences, Corvinus University of Budapest, Hungary \\ \email{peter.biro@krtk.mta.hu} \and Radu-Stefan
 Mincu \and Alexandru Popa \at Department of Computer Science, University of Bucharest, Romania \and Utkarsh Verma \at
Department of Industrial Engineering and Operations Research, IIT Bombay, India\\
}
\maketitle
\begin{abstract}
In kidney exchange programmes patients with end-stage renal failure may exchange their willing, but incompatible living donors among each other. National kidney exchange programmes are in operation in ten European countries, and some of them have already conducted international exchanges through regulated collaborations. The exchanges are selected by conducting regular matching runs (typically every three months) according to well-defined constraints and optimisation criteria, which may differ across countries. In this work we give integer programming formulations for solving international kidney exchange problems, where the optimisation goals and constraints may be different in the participating countries and various feasibility criteria may apply for the international cycles and chains. We also conduct simulations showing the long-run effects of international collaborations for different pools and under various national restrictions and objectives.
\end{abstract}
\section{Introduction}

When an end-stage kidney patient has a willing, but incompatible living donor, then in many countries this patient can exchange his/her donor for a compatible one in a so-called kidney exchange programme (KEP). The first national kidney exchange programme was established in 2004 in the Netherlands in Europe \cite{de_etal2005}. Currently there are ten countries with operating programmes in Europe \cite{Biro_etal2018}, the largest being the UK programme \cite{MO14}.

Typically the matching runs are conducted every three months on pools with around 50-300 patient-donor pairs. The so-called virtual compatibility graph represents the patient-donor pairs with nodes and an arc represents a possible donation between the corresponding donor and patient, that is found compatible in a virtual crossmatch test. The exchange cycles are selected by well-defined optimisation rules, that can vary across countries. The most important constraints are the upper limits on the length of exchange cycles, for examples, two in France, three in the UK and Spain, and four in the Netherlands \cite{Biro_etal2018}. The main goal of the optimisation in Europe is to facilitate as many transplants as possible, i.e.\ to maximise the number of nodes covered in the compatibility graph by independent cycles. The corresponding computational problem for cycle-length limits three or more is NP-hard, and the standard solution technique used is integer programming \cite{ABS07}.

International kidney exchanges have already started in Europe between Austria and Czech Republic \cite{Bo_etal17} since 2016, between Portugal, Spain and Italy since summer 2018, and between Sweden, Norway and Denmark in the Scandiatransplant programme (STEP), built on the Swedish initiative \cite{AK16}. The Italy-Portugal-Spain collaboration is organised in a consecutive fashion, first the national runs are conducted and then the international exchanges are sought for the remaining patient-donor pairs. A related game-theoretical model has been studied in \cite{Ca_etal2017}. In the Scandinavian programme and in the Prague-Wien collaboration, the protocol is to find an overall optimum for the joint pool. In the latter situation, the fairness of the solution for the countries involved can be seen as an important requirement, which was studied in \cite{KNPV} with extensive long-term simulations by proposing the usage of a compensation scheme among the countries.

A similar situation arises in the US kidney exchange problem, where the transplant centres are the strategic agents \cite{AR12,AR14,As_etal15,TP15}. Most of the scientific studies focus on the problem of incentivising the hospitals to register all of their patient-donor pairs, and not only the hard-to-match ones, which often happens in practice as the majority of transplants are still conducted within transplant centres, outside of the three national schemes \cite{Ag_etal18}. Novel suggestions of credit based schemes have also been studied \cite{Ha_etal15,Ag_etal18}, and a similar system has been implemented in the National Kidney Registry, which is the largest in the volume of transplants among the three nationwide kidney exchange programme in the US.

In this study we focus on the collaboration of countries and a key aspect of this collaboration is the assumption that they all follow commonly agreed protocols. As such, there is no need to incentivise countries to register their patient-donor pairs, unlike for the American hospitals. We will however compare the consecutive and the joint pool scenarios in our simulations, as these are both used in practice. We will not consider compensations, or any strategic issues, but we will allow the countries to have different constraints and goals with regard to the cycles and chains they may be involved in. In particular, we will compare the benefits of the countries from international collaborations when they have different upper bounds on their national cycles, and thus also possible different constraints on the segments of the international cycles they are participating in. As an example we mention the Austro-Czech cooperation, where Austria requires to have all exchanges simultaneously, so they allow short national cycles and short segments only, whilst in Czech Republic longer non-simultaneous cycles and chains are also allowed. We formulate novel IP models for dealing with potentially diverse constraints and goals in international kidney exchange programmes and we test two-country cooperation scenarios under different assumptions over their constraints, the possibility of having chains triggered by altruistic donors, and the sizes of their pools.

\section{Model of international kidney exchanges}

In a standard kidney exchange problem, we are given a directed \emph{compatibility graph} $D(V,A)$, where the nodes $V=\{1,2,\dots n\}$ correspond to patient-donor pairs and there is an arc $(i,j)$ if the donor of pair $i$ is compatible with the patient of pair $j$. Furthermore we have a non-negative weight function $w$ on the arcs, where $w_{i,j}$ denotes the weight of arc $(i,j)$, representing the value of the transplantation. In most application the primary concern is to save as many patients as possible, so the weight of each arc is simply equal to one. We will also focus on this case in our simulations.

Let $\mathcal{C}$ denote the set of cycles allowed in $D$, which are typically allowed to be of length at most $K$. The solution of a classical kidney exchange problem is a set of disjoint cycles of $\mathcal{C}$, i.e., a cycle-packing in $D$. For cycle $c\in\mathcal{C}$, let $A(c)$ denote the set of arcs in $c$ and $V(c)$ denote the set of nodes covered by $c$.\footnote{In addition, we can also consider \emph{altruistic donors}, in which case we separate the node set into patient-donor pairs $V_p$ and altruistic donors $V_a$, so $V=V_p\cup V_a$. The solution may contain exchange cycles and chains triggered by altruistic donors. The latter can be conducted non-simultaneously, so different restrictions may apply for them. In this paper we focus on cycles, but we note that one can reduce the problem of finding chains to the problem of finding cycles by adding artificial patients to the altruistic donors, who are compatible with all donors.}

In an international kidney exchange programme multiple countries ($N$) are involved in the exchange, so the set of nodes is partitioned into $V=V^1\cup V^2\cup\dots\cup V^N$, where $V^k$ is the set of patient-donor pairs in country $k$. We have the following modifications of the classical problem. Let $A^k$ denote the arcs pointing to $V^k$ (so the donations to patients in country $k$). Note that $A=A^1\cup A^2\cup\dots\cup A^N$. The weights of the arcs in $A^k$ should reflect the preferences of country $k$. (We may assume that these are scaled, e.g., by having the same average scores in order not to bias the overall optimal solution towards some countries.) Finally, let $A^{\mathcal{N}}$ and $A^{\mathcal{I}}$ denote the national and international donations, i.e.,\ $A=A^{\mathcal{N}}\cup A^{\mathcal{I}}$.

In a global optimal solution, small cycles within the countries can have different requirement than international cycles. Therefore we separate the two sets of cycles into $\mathcal{C}=\mathcal{C^N}\cup\mathcal{C^I}$, where $\mathcal{C^N}$ is the set of national cycles and $\mathcal{C^I}$ is the set of international cycles. We call the national parts of an international cycle \emph{segments}, where a segment is a path within a country, and the segments are linked by international arcs. A $l$-segment is a path of length $l-1$, with all the $l$ nodes belonging to the same country. Let $\mathcal{S}$ denote the set of all possible segments, and let $\mathcal{S}^k$ denote the set of segments allowed in country $k$. For $s\in\mathcal{S}$, let $A(s)$ denote the set of (national) arcs and let $V(s)$ denote the set of nodes covered (in the same country). Note that a segment may also consist of a single node, which corresponds to the case when an international donation is immediately followed by another international donation.

We can have the following natural restrictions on the national and international cycles.\footnote{In addition we can also have different constraints for altruistic chains, and we may require that an international chain may have to end in the same country where it started.} We may have upper limits on the

\begin{enumerate}
\item[C1)] total length of an international cycle
\item[C2)] length of national cycles for each country (possibly different)
\item[C3)] length of segments in international cycles for each country (possibly different)
\item[C4)] number of segments in a country within an international cycle
\item[C5)] number of patient-donor pairs from a country in an international cycle
\item[C6)] number of countries involved in an international cycle
\end{enumerate}

\section{Integer programming formulations}
\label{sec:integerform}

First we describe the two classical IP formulations of Abraham et al. \cite{ABS07}, the edge-formulation and the cycle-formulation. We will build our general IP solutions on these.

\subsection{Basic edge-formulation}

We introduce a binary variable $y_{i,j}$ for each arc $(i,j)$. Finding a maximum (value) solution with cycles of length at most $K$ can be

\begin{equation}
\label{obj:edge}
max \sum_{i,j}w_{i,j}y_{i,j}
\end{equation}

such that

\begin{equation}
\label{eq:Kickhoff}
\sum_{i:(i,j)\in A}y_{i,j}=\sum_{k:(j,k)\in A}y_{j,k}  \hspace{2em} \forall j\in V
\end{equation}

\begin{equation}
\label{eq:path}
\sum_{j:(i,j)\in A}y_{i,j}\leq 1  \hspace{2em} \forall i\in V
\end{equation}

\begin{equation}
\label{eq:length}
\sum_{(i,j)\in A(p)}y_{i,j}\leq K-1  \hspace{2em} \forall p\in\mathcal{P}^K
\end{equation}

where $\mathcal{P}^K$ denotes the set of $K$-length proper directed paths (i.e.\ which are not cycles), and $A(p)$ denotes the set of arcs in $p$.

\subsection{Basic cycle-formulation}

We introduce a binary variable $x_c$ for each cycle $c\in\mathcal{C}$. The weight of a cycle $c$ is denoted by $w_c$, which can be taken as the sum of the edge-weights in the cycle, or can be defined differently.

\begin{equation}
\label{obj:cycle}
\max \sum_{c\in\mathcal{C}} w_c x_c
\end{equation}

such that

\begin{equation}
\label{eq:cycle-feasible}
\sum_{c: i\in V(c), c\in\mathcal{C}}x_c\leq 1\  \hspace{2em} \forall i\in V
\end{equation}

One can use the cycle-formulation in our international setting after carefully searching and selecting the potential national and international cycles. We used this technique to solve the two-country problems, described in the last section. \textbf{We defer the description of the cycle-search algorithm to the Appendix}.

\subsection{Linking the cycle-, segment-, and edge-variables}

When the pre-selection of international cycles is complicated then we can use a mixed formulation, where the national cycles are represented by binary variables, and the international cycles are decomposed into allowable segments. We show how to link the cycle and segment variable with the edge variables, and thus enforce various constraints for different countries. Let $z_s$ be a binary variable of segment $s\in\mathcal{S}$. Suppose that we only work with national cycles $\mathcal{C^N}$ with binary variables for all of them, but we do not have variable for international cycles.

Besides the basic feasibility cycle-constraints \eqref{eq:cycle-feasible} and edge-constraints \eqref{eq:Kickhoff}, we need the following sets of equations.

\begin{equation}
\sum_{j:(i,j)\in A} y_{i,j} \leq \sum_{c:i\in V(c)} x_c +\sum_{s:i\in V(s)} z_s\leq 1, \forall i \in V \label{eq:cycle-edge3}
\end{equation}

The above condition \eqref{eq:cycle-edge3} enforces that we can only cover a node by either a national cycle or by a segment.

\begin{equation}
|A(c)|\cdot x_c\leq \sum_{(i,j)\in A(c)}y_{i,j} \hspace{2em} \forall c\in \mathcal{C^N} \label{eq:cycle-edge1}
\end{equation}

\begin{equation}
\sum_{(i,j)\in A(c)}y_{i,j}-|A(c)|+1\leq x_c \hspace{2em} \forall c\in \mathcal{C^N} \label{eq:cycle-edge2}
\end{equation}

Conditions \eqref{eq:cycle-edge1} and \eqref{eq:cycle-edge2} imply the inclusion of all the edge-variables in a national cycle if and only if that cycle is selected in the solution.

Let $e^{+}(i)=\sum_{(j,i)\in A^{\mathcal{I}}}y_{j,i}$ and let $e^{-}(i)=\sum_{(i,j)\in A^{\mathcal{I}}}y_{i,j}$ two new indicator variables showing whether node $i$ is receiving or giving a kidney in an international exchange. We set the following condition for each segment $s$ with starting node $u$ and terminal node $v$.
\begin{equation}
(|A(s)|+2)\cdot z_s\leq \sum_{(i,j)\in A(s)}y_{i,j}+e^{+}(u)+e^{-}(v)\hspace{2em} \forall s\in \mathcal{S} \label{eq:segment-edge1}
\end{equation}
\begin{equation}
\sum_{(i,j)\in A(s)}y_{i,j}+e^{+}(u)+e^{-}(v) - |A(s)|-1 \leq z_s \hspace{2em} \forall s\in \mathcal{S} \label{eq:segment-edge2}
\end{equation}
Conditions \eqref{eq:segment-edge1} and \eqref{eq:segment-edge2} imply the inclusion of all the edge-variables of a segment if and only if that segment is selected as part of an international cycle in the solution.

\subsection{Satisfying the special constraints}

Here we give possible solutions for enforcing the requirements with IP formulations. As we noted, by a careful cycle-search algorithm we can always satisfy all of the conditions as long as the short upper bounds on the lengths of cycles exists. However, if a country allows unbounded length cycles or chains then we shall use edge-variables as well, combined with cycle and segment-variables. Furthermore, the usage of edge- and segment-variables can also be useful to simplify the cycle-search algorithms and rule out the infeasible cycles by edge-formulation constraints instead. Therefore, in the following description we explain how to use the variants of the previously provided edge-formulations for achieving the required conditions.

\subsubsection*{C1) total length of an international cycle}

In case we have only this requirement may use the basic edge-formulation to satisfy this condition. However, if we also have restrictions on the national cycles then we need an extended model that we will describe in the following point.

\subsubsection*{C2) length of national cycles for each country}

Regarding the edge-formulations, we may already need a more sophisticated formula, in case the international cycle would allow a longer segment in a country than the maximum length of the national cycles in that country. For this purpose, we can introduce two new edge variables for each arc $(i,j)$, $\hat{y}_{i,j}$ and $\check{y}_{i,j}$, where the former denotes whether this edge is used in a national cycle and the latter denotes whether it is used in an international cycle. We have $y_{i,j}=\hat{y}_{i,j}+\check{y}_{i,j}$, and for any arc in $A^{\mathcal{I}}$ we do not have $\hat{y}_{i,j}$, as this arc cannot be involved in a national cycle. We also have to ensure consistency, so condition \eqref{eq:Kickhoff} should be written up separately for national and international edge-variables for each node as follows.

\begin{equation}
\label{eq:Kickhoff_nat}
\sum_{i\in V}\hat{y}_{i,j}=\sum_{k\in V}\hat{y}_{j,k}  \hspace{2em} \forall j\in V
\end{equation}

\begin{equation}
\label{eq:Kickhoff_int}
\sum_{i\in V}\check{y}_{i,j}=\sum_{k\in V}\check{y}_{j,k}  \hspace{2em} \forall j\in V
\end{equation}

The reason behind this separation is to allow different upper bounds on the lengths of the national and international cycles. For instance, if the international cycles have upper bound $K$ and in the meantime country $k$ have upper bound $K^k$ on the length of its national cycles then we can achieve both by the following pair of conditions.

\begin{equation}
\label{eq:length_int}
\sum_{(i,j)\in A(p)}\check{y}_{i,j}\leq K-1  \hspace{2em} \forall p\in\mathcal{P}^K \mbox{ in } D
\end{equation}

\begin{equation}
\label{eq:length_int}
\sum_{(i,j)\in A(p)}\hat{y}_{i,j}\leq K^k-1  \hspace{2em} \forall p\in\mathcal{P}^{K^k} \mbox{ in } D^k
\end{equation}

It is also important to note that after this separation we shall set the constraints used for making the connections between the edge-, cycle- and segment-variables accordingly. Namely, in the constraints for national cycles, \eqref{eq:cycle-edge1} and \eqref{eq:cycle-edge2}, we shall replace $y_{i,j}$ with $\hat{y}_{i,j}$, whilst in the constraints for segments, \eqref{eq:segment-edge1} and \eqref{eq:segment-edge2}, we shall replace $y_{i,j}$ with $\check{y}_{i,j}$.

\subsubsection*{C3) length of segments in international cycles for each country}

To bound the length of national segments in international cycles, we can either search the allowable segments and define segment variables for them, or we can satisfy these constraints by using the international edge-variables. For instance if $L^k$ is the upper bound for the length of segments in country $k$ then we need the following modified edge-constraint.

\begin{equation}
\label{eq:length_int}
\sum_{(i,j)\in A(p)}\check{y}_{i,j}\leq L^k  \hspace{2em} \forall p\in\mathcal{P}^{L^k} \mbox{ in } D^k
\end{equation}

\subsubsection*{C4) number of segments in a country within an international cycle}

Following the idea of \cite{ConstantinoKVR13} used for providing a compact formulation for the basic problem, we propose to define layers in order to separate the international cycles from each other. This will facilitate a simple way to set restrictions for this and the final two points.

Suppose that we can have at most $T$ international cycles in the solution (e.g., $T\leq |V|/2$). For every $t\in [1\dots T]$ we define binary edge-variables $y_{i,j}^t$ denoting whether that edge is included in the $t$-th cycle. We set $\check{y}_{i,j}=\sum_{t\in [1\dots T]}y_{i,j}^t$. We replicate \eqref{eq:Kickhoff} for each layer, as follows.

\begin{equation}
\label{eq:Kickhoff_t}
\sum_{i\in V}y_{i,j}^t=\sum_{k\in V}y_{j,k}^t  \hspace{2em} \forall j\in V, t\in [1\dots T]
\end{equation}

Now, we can restrict the number of segments from the same country $k$ to be less than or equal to $\lambda^k$ with the following conditions.

\begin{equation}
\label{eq:country_segments}
\sum_{(i,j)\in A^{\mathcal{I}}\cap A^{k}}y_{i,j}^t\leq \lambda^k  \hspace{2em} \forall t\in [1\dots T]
\end{equation}

Note that the above formulation does not rule out the possibility of having multiple international cycles in one layer, however, neither of them can have more than $\lambda^k$ segments in country $k$. Furthermore, if we want to enforce that an international cycle can only have one segment from each country in a two-country programme then we can simplify the constraints, as follows. We can separate the layers according to the nodes from the first country that has outgoing arcs to the second country. Let the $t$-th such node, $i$, have only variables $y^t_{i,j}$ for outgoing international transplants, which means that this node can only be involved in an international cycle at the $t$-th layer. See more about this idea at the end of this Section, where we describe the implementation of the bounder-unbounded two-country case (i.e., Austria-Czech Republic collaboration) in details.

\subsubsection*{C5) number of patient-donor pairs from a country in an international cycle}

By using again the above defined layers, we can easily set restrictions on the number of pairs involved in an international cycles from one country. If this upper bound is $\beta^k$ for country $k$ then we can enforce this condition with the following constraints.

\begin{equation}
\label{eq:country_pairs}
\sum_{(i,j)\in A^{k}}y_{i,j}^t\leq \beta^k  \hspace{2em} \forall t\in [1\dots T]
\end{equation}

\subsubsection*{C6) number of countries involved in an international cycle}

To achieve this restriction with edge-variables, we shall define a new indicator variable for each country-layer pair, which shows whether this country is involved in an international cycle in that layer. Thus, let us define a binary variable $b_k^t$ for each country $k$ and layer $t$, which must be one if this country is involved in an international cycle at this layer. This can be achieved with the following constraints.

\begin{equation}
\label{eq:country_number}
\sum_{(i,j)\in A^{\mathcal{I}}\cap A^{k}}y_{i,j}^t\leq b_k^t\cdot |V^k|  \hspace{2em} \forall k\in [1\dots n], t\in [1\dots T]
\end{equation}

This implies that if there is an ingoing international transplant to country $k$ at layer $t$ then $b_k^t$ must be equal to one (but otherwise it can be zero). If the number of countries allowed in an international cycle is $\gamma$ then we can enforce this condition with the following constraint.

\begin{equation}
\label{eq:country_number_bound}
\sum_{k}b_k^t\leq \gamma  \hspace{2em} \forall t\in [1\dots T]
\end{equation}

\subsubsection*{Adding altruistic chains}

We already noted that if we add the possibility of having altruistic chains, where their lengths is similarly bounded as the cycles then we can simply treat them as cycles by adding a dummy arc to the imaginary patient of the altruistic donor from all donors. However, the length restrictions are different for the altruistic chains (which can be reasonable, since these can be conducted non-simultaneously), in fact, they can even be considered unbounded (aka.\ \emph{never ending chains} \cite{RKP09}). In this case we shall introduce new edge-variables corresponding to the possibility of conducting a transplant in an altruistic chain. Furthermore, we may even separate them to variables for national and international chains, if they have different restrictions.

\subsection*{An example: the case of Austria and Czech Republic}\label{sec:A-CZ}

As an example, we describe the full IP-model for a problem setting that is representing the cooperation between Austria and Czech Republic. Here Austria allows only short national cycles and short segments, whilst Czech Republic allows unbounded national cycles and unbounded segments in international cycles, and we require that every international cycle has only one segment in each country. This IP formulation was implemented and used in the simulations for the two-country cases, that we described in Section \ref{sec:simulation}.

Let $V^1$ denote the first country with bounded length cycles and segments (Austria) and let $V^2$ denote the second country with unbounded length national cycles and segments (Czech Republic). We introduce cycle and segment variables for $V^1$ (but not for $V^2$) and we introduce edge-variables for $A^2$ and $A^{\mathcal{I}}$ (but not for the internal edges in $V^1$). Furthermore, we introduce layers $1\dots T$, such that $T$ is the number of nodes in $V^1$ that has outgoing edge to $V^2$. For all the edges with edge variables $y_{i,j}$ we introduce also the layer variables $y^t_{i,j}$ for every $t\in [1\dots T]$, except the edges from $V^1$ to $V^2$. For these edges we introduce only one layer variable, as follows. Let $i$ be the $t(i)$-th node in $V^1$ that has outgoing edge to $V^2$. We introduce a binary variable $y^{t(i)}_{i,j}$ for every arc $(i,j)\in A^{\mathcal{I}}$, but these edges will not have other layer variables. This will ensure that every international cycle has only one segment in each country. Below we describe the complete IP formulation for this special setting.

\begin{subequations}
\begin{align}
\max  \sum\limits_{c \in \mathcal{C}} |c|x_c + \sum\limits_{(i,j) \in A}  y_{i,j} & + \sum\limits_{s \in \mathcal{S}} |s|z_{s}   \nonumber\\
\text{ subject to:} \sum\limits_{t\in [1\dots T]} y_{i,j}^{t} & = y_{i,j} \hspace{2em} \forall (i,j)\in A\\
\sum\limits_{i\in V} y_{i,j}^{t} & = \sum\limits_{k\in V} y_{j,k}^{t} \hspace{2em} \forall t\in [1\dots T], j\in V \\
\sum\limits_{j\in V} y_{i,j} & \leq 1 \hspace{2em} \forall i\in V \\
\sum\limits_{c: i\in V(c)} x_c  + \sum\limits_{s: i\in V(s)} z_{s} & \leq 1 \hspace{2em} \forall i\in V^1 \\
\sum\limits_{k\in V^2} y_{k,i} & = \sum\limits_{s: s = (i,\dots)} z_{s} \hspace{2em} \forall i\in V^1\\
\sum\limits_{k\in V^2} y_{j,k} & = \sum\limits_{s: s = (\dots, j)} z_{s} \hspace{2em} \forall i\in V^1\\
\sum\limits_{k\in V^2} y_{k,i}^{t(j)} + \sum\limits_{l\in V^2} y_{j,l}^{t(j)} & \geq  2 z_s \hspace{2em} \forall i,j\in V^1, s=(i, \dots, j)
\end{align}
\end{subequations}

\section{Simulations}
\label{sec:simulation}
We conduct long-term simulations with agents arriving and leaving the pool such as in \cite{San17}, (e.g. with 3-months matching runs for 3
years). We restrict our attention to maximising the size of the solutions, and do not consider scoring methods or any other objective. We conduct the simulations for the two-country case, as the effects of the cooperation for a country can already be tested on this simple setting.

\subsubsection*{Cooperation policies:}

We consider three basic policies for international cooperation:

\begin{enumerate}[a)]
\item \emph{no cooperation}: each country conducts its matching run separately.
\item \emph{consecutive matching}: each country conducts its national run first and then the
international run is conducted for the remaining pairs (as done in the
Spanish-Italian-Portuguese programme).
\item \emph{merged pools}: where the countries register all their pairs in a merged pool and a global optimisation is conducted (e.g. Austria-Czech Republic and Sweden-Norway-Denmark) while respecting local policy restrictions.
\end{enumerate}

\begin{figure}[b]
\caption{Graphical representation of solutions for the first KEP stage in one of the instances: altruist donors are at the top, patient-donor pairs form circles for each country and arcs represent transplants. Left side, individual KEPs: 13/16 patients receive transplants in the small country, 28/38 patients in the large country are transplanted. Right side, merged KEP: the numbers are 15/16 for the small country, and 32/38 for the large one. \vspace{2mm}} \label{fig:example}
\includegraphics[width=0.45\textwidth]{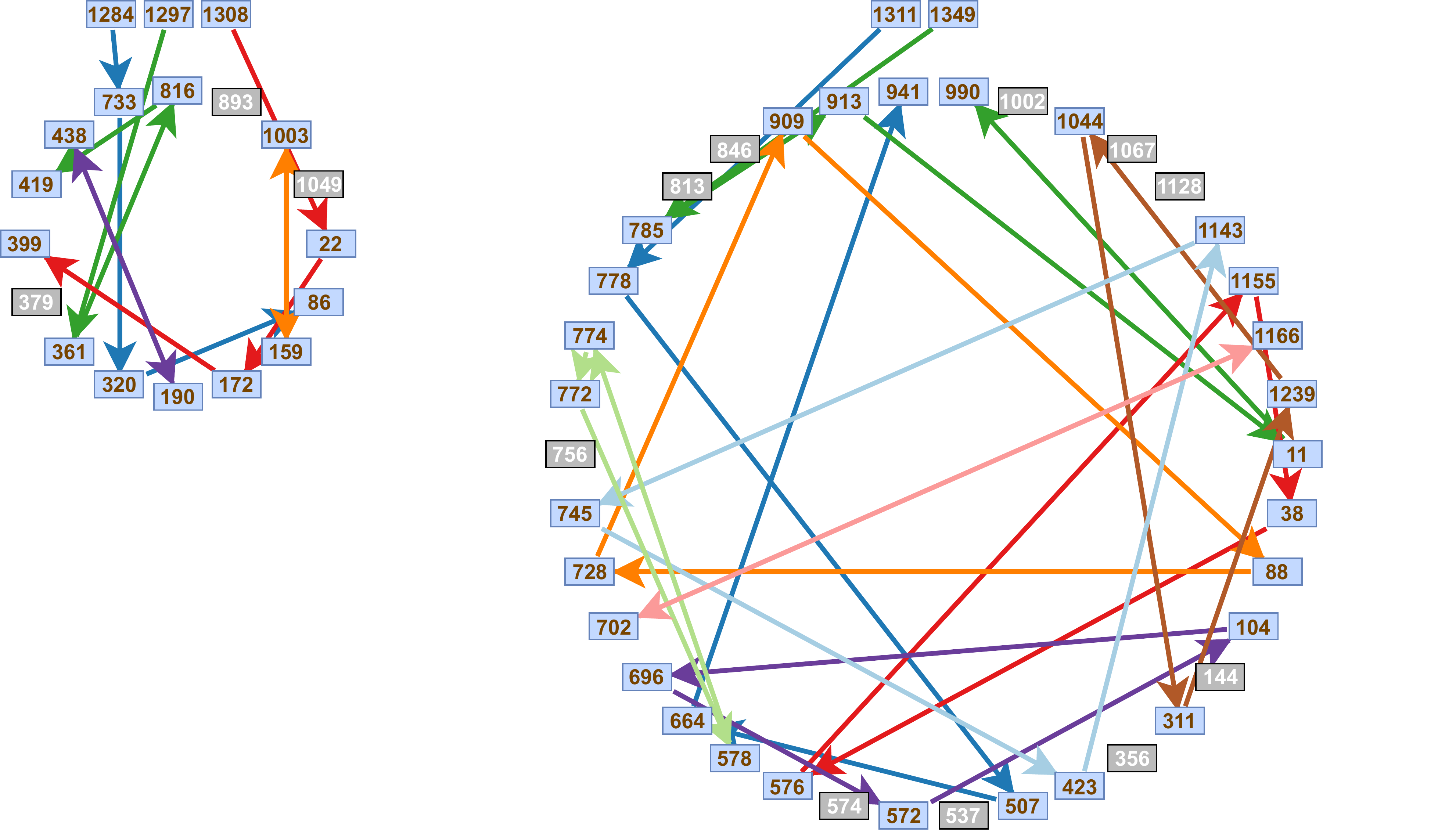}
\includegraphics[width=0.55\textwidth]{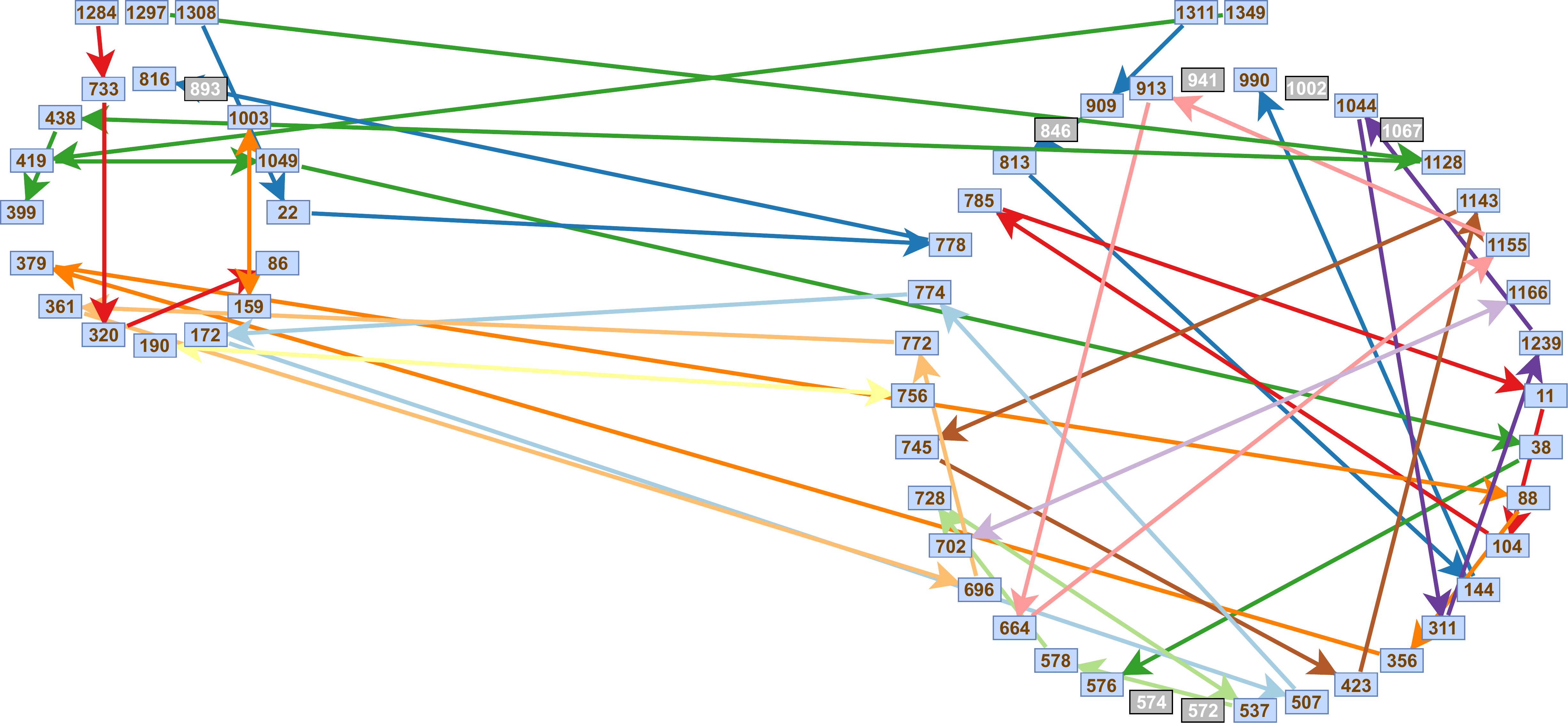}
\vspace{-8mm}
\end{figure}

\subsubsection*{Differences between countries}

Our main interest is to evaluate how the differences between the national pools and the optimisation constraints effect the benefits of the countries involved in the collaboration.

Patient availability constraints:
\begin{enumerate}
  \item[1.] size of pools: size ratios being 1:1, 3:4, 1:2, 1:4. 
    \end{enumerate}

Local policy constraints:
    \begin{enumerate}
  \item[2.] different upper bounds on the length of the national cycles: 2:2, 2:3, 2:4, 2:$\infty$, 3:3, 3:4, 3:$\infty$.
    \end{enumerate}

International policy constraints:
    \begin{enumerate}
  \item[3.] the size of the international cycles: it is set to be the larger individual cycle bound.
  \item[4.] the length of the segments in each country: it is set to one less than the individual cycle bound
  \item[5.] the number of segments in each country: it is unrestricted with the exception of the case of collaboration with a country with unbounded cycles where we consider only one segment for each country per international cycle.

\end{enumerate}

In the summary tables we show how much each country is benefiting from the cooperation, as compared to the baseline setting (no cooperation).

\section{Case Study: Simulations involving two countries}
To determine the benefits of international KEPs we conduct a case study involving two countries which aim to develop a joint KEP and are concerned about the advantages and disadvantages of cooperation between their KEPs. We compare the individual benefits from the no cooperation case to the consecutive matching and merged pool scenarios.

The simulation involves 20 instances each containing the compatibility information for 1000 patient-donor pairs. For the sake of simplicity, although our model can handle altruist donors (who start chains), for the case study we only consider the cyclic exchanges among incompatible patient-donor pairs. The reason for this choice is to be able to more clearly see the difference in cooperation (because adding chains may dramatically increase the number of transplants in any given KEP stage). The length of the considered time-frame for the simulated kidney exchange scenarios is 3 years with matching runs scheduled every 3 months for each instance. Every agent is assigned an uniformly distributed arrival time to the KEP and the patient-donor pairs stay in the KEP for a maximum of 1 year (or 4 matching runs) after which they leave the programme (which means that they opt for an alternative solution, such as having a direct ABO-incompatible transplant after desensitisation or getting a deceased organ).

To understand the importance of the parameters of collaboration presented in Section \ref{sec:simulation} we conduct a sensitivity analysis by considering most of the combinations possible. We illustrate the impact of both programme pool size and local policy constraints on the cooperation of KEPs by:
\begin{enumerate}
\item examining the stage-by-stage amount of transplants and drop-out patients of the programme in each country (see Figure \ref{fig:perstageplots}).
\item examining the total number of transplants resulting after 3 years for each participant country (see Figures \ref{fig:surf0:samepool}-\ref{fig:surf4:bound4}).
\end{enumerate}

\subsubsection*{Explanation of the figures and tables}

In the first row of Figure \ref{fig:perstageplots} we consider a baseline scenario where two countries have the same pool size (on average 41.6 patients arrive per stage or 500 during 3 years) and same constraints (local and international cycle bound 3, no altruist donors).

Then, in the second row of Figure \ref{fig:perstageplots} we demonstrate the effects of local KEP constraints on the matched pairs by only changing the local constraint in country 2 (now has local cycle bound 2 and international cooperation can be done with 2-cycles and 3-cycles with only one participant patient-donor pair in country 2). Note that the patient compatibility and arrival time to the programme was not modified.


In the third row of Figure \ref{fig:perstageplots}, we show the impact of pool size on the number of transplants by only removing half of the patients inside the pool of country 2 (they never register to the KEP). Note that the patient compatibility and arrival time to the programme for the remaining patients was not modified.

Finally, the bottom row of Figure \ref{fig:perstageplots} shows the relative improvement in the number of transplants defined as: the number of transplants obtainable by merged KEP divided by the number of transplants that would have been achieved by the respective country's local KEP.


In every scenario the objective is simply to maximise the number of transplants. There are three settings for collaboration: no cooperation (i.e. separate KEPs, baseline scenario), consecutive matchings (each country runs a local KEP optimisation and then the remaining patient pools enter a joint KEP) and full collaboration (a single KEP for both countries).

Figures \ref{fig:surf0:samepool}-\ref{fig:surf4:bound4} describe the improvement from merged KEP collaboration in the number of transplants over the number of transplants each country can achieve by themselves (we define this value as the \textit{benefit} of collaboration). These figures are a result of a sensitivity analysis experiment where the effects of local country bounds and size limitations are explored. The information is also displayed in numeric format in Tables \ref{table1:bound} and \ref{table2:sizebound}. The listed values represent the number of performed transplants (for each country) measured at the end of the 3 year programme and the values are averages of 20 instances (each instance runs for 12 KEP stages).

\subsubsection*{Experimental results discussion}

\begin{figure}[H]
\caption{Country 1 has cycle bound 2, 100\% pool size}
\label{fig:surf2:bound2}
\includegraphics[width=0.495\textwidth]{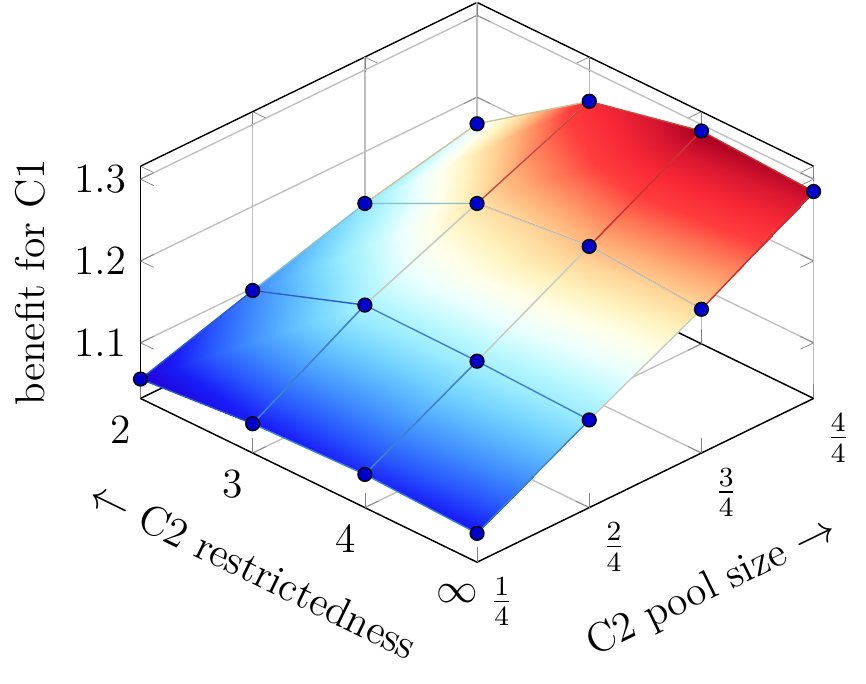}
\includegraphics[width=0.495\textwidth]{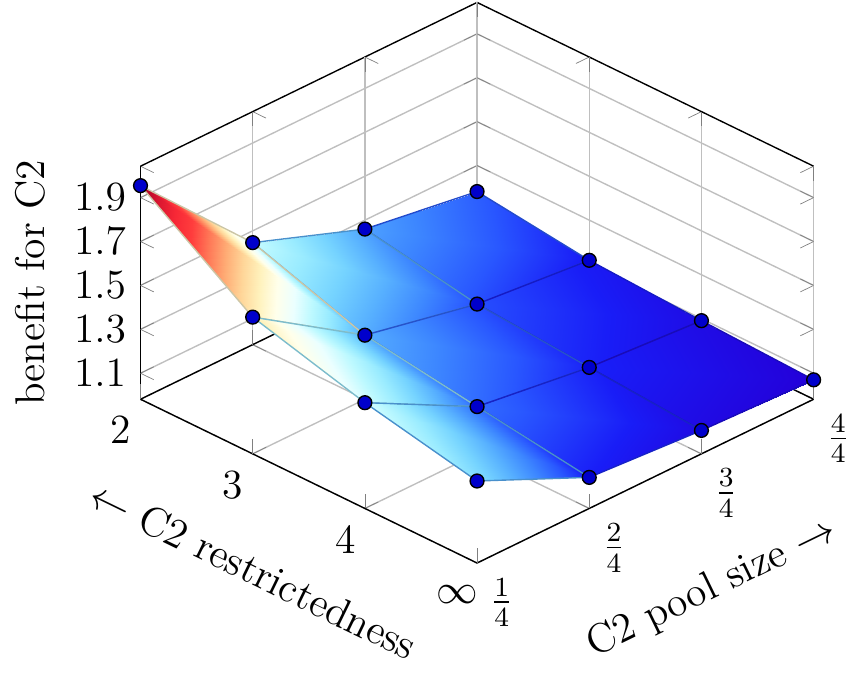}
\end{figure}

\begin{figure}[H]
\caption{Country 1 has cycle bound 3, 100\% pool size}
\label{fig:surf3:bound3}
\includegraphics[width=0.495\textwidth]{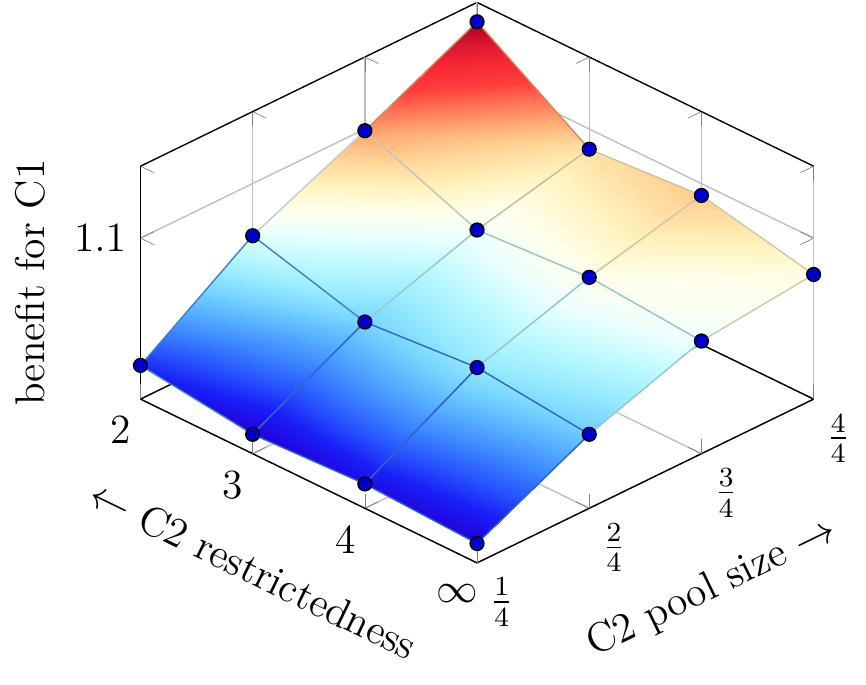}
\includegraphics[width=0.495\textwidth]{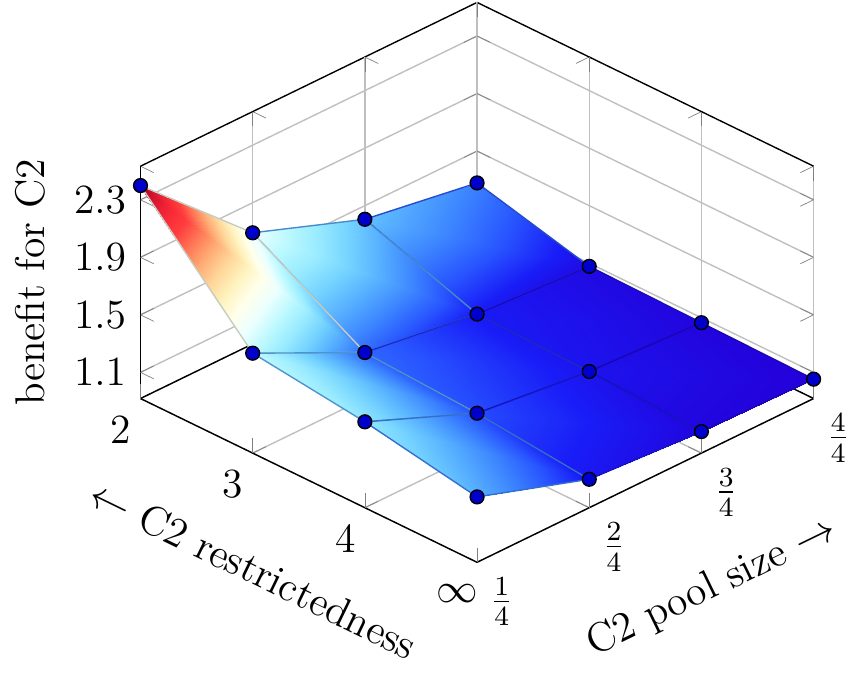}
\end{figure}

\begin{figure}[H]
\caption{Country 1 has cycle bound 4, 100\% pool size}
\label{fig:surf4:bound4}
\includegraphics[width=0.495\textwidth]{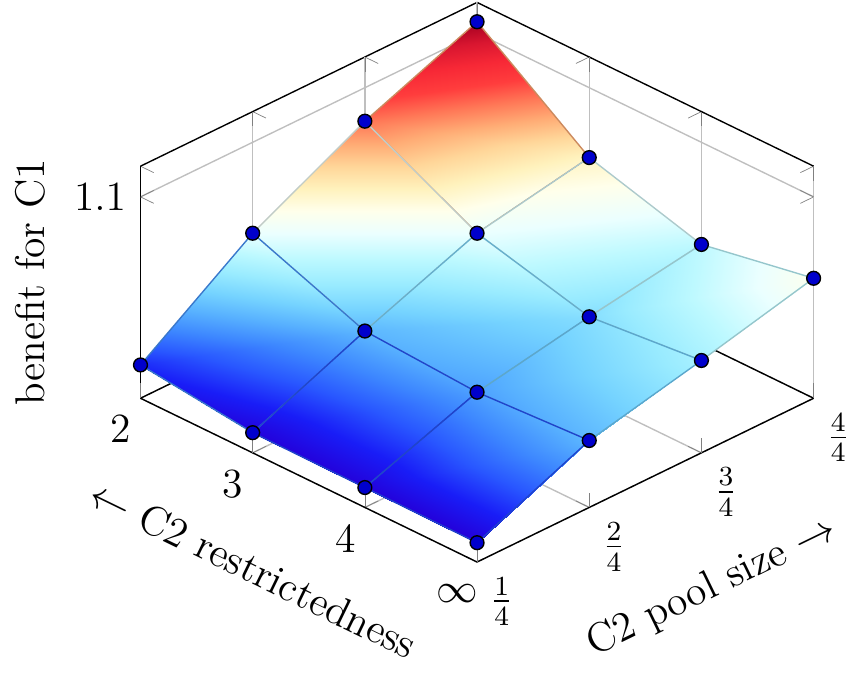}
\includegraphics[width=0.495\textwidth]{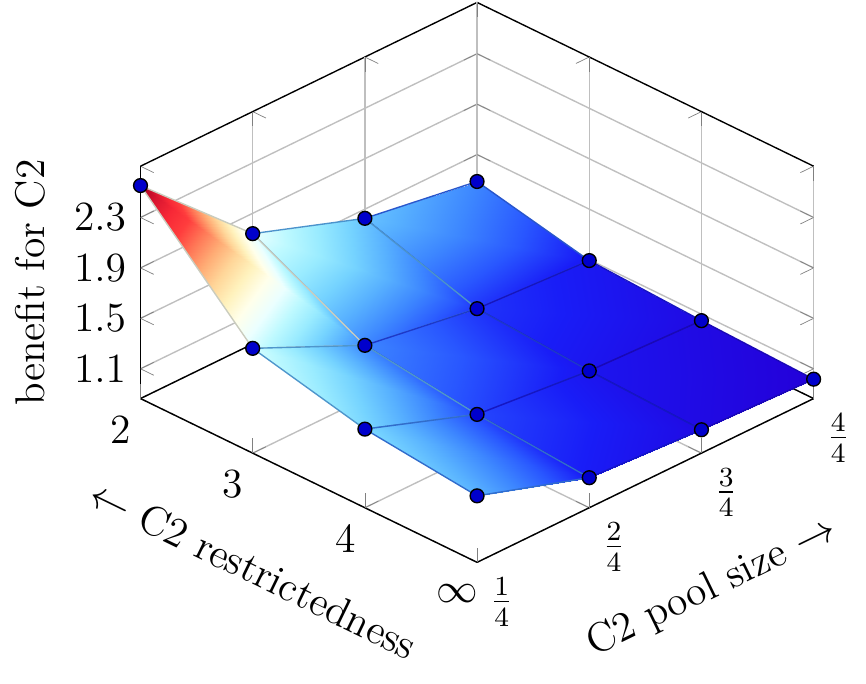}
\end{figure}

\begin{figure}[H]
\caption{Plots detailing the amount of transplants (solid line) and dropouts (dotted line) during each stage of the KEP and separately for each country (country 1 on the left side and country 2 on the right). The first three rows each containin two plots showcasing a different scenario: (1:1 size, 3:3 bounds), (1:1 size, 3:2 bounds) and (2:1 size, 3:3 bounds), respectively. The bottom row describes individual relative improvement $\frac{\text{merged transplants}}{\text{local transplants}}$.} \label{fig:perstageplots}
\includegraphics[width=0.495\textwidth]{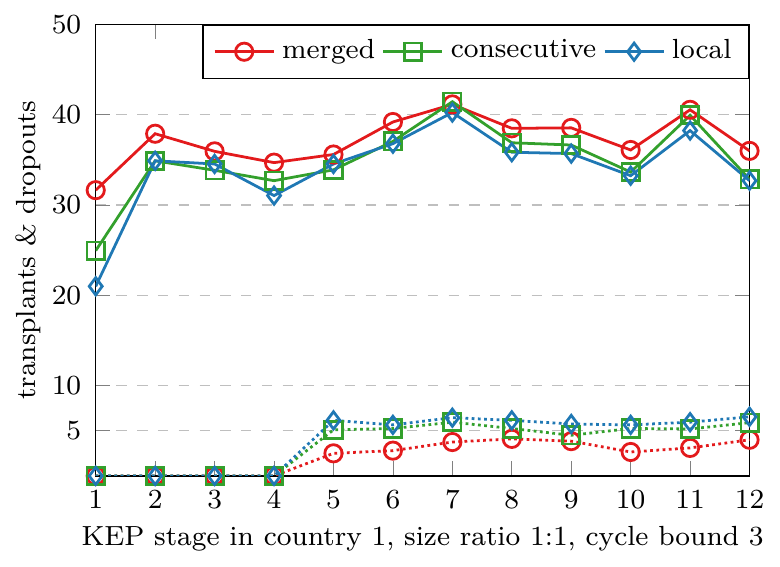}
\includegraphics[width=0.495\textwidth]{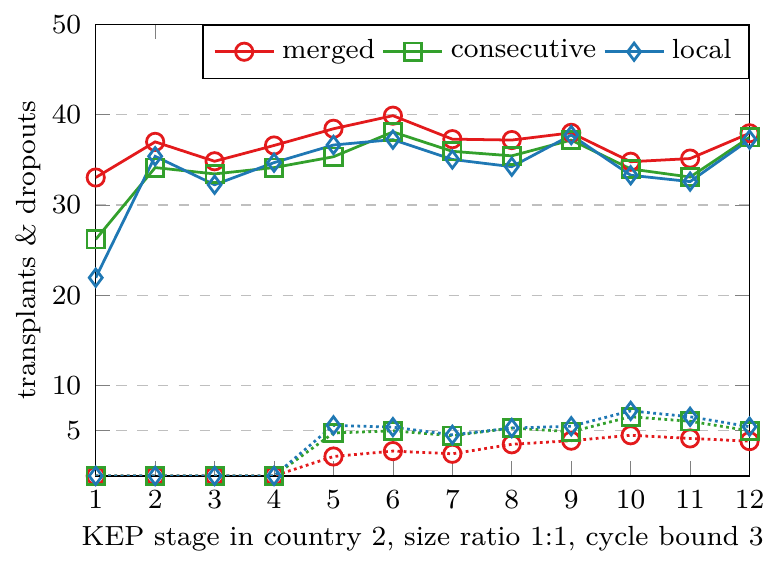}

\includegraphics[width=0.495\textwidth]{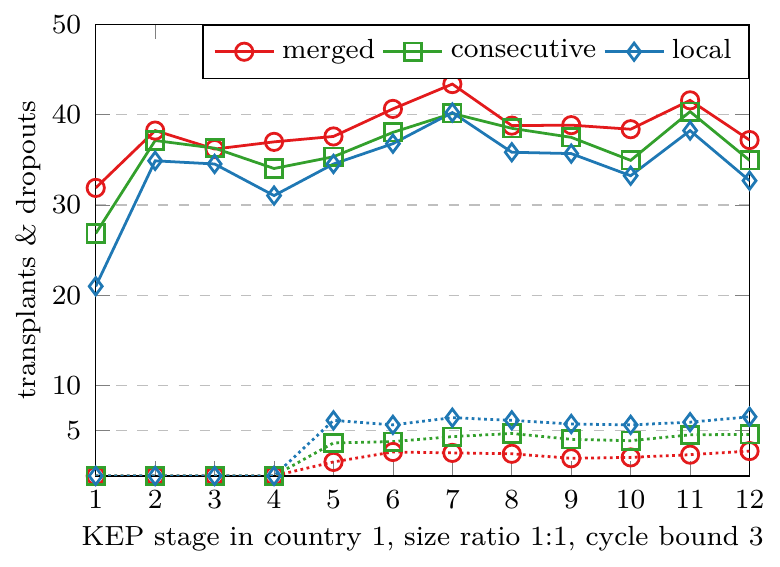}
\includegraphics[width=0.495\textwidth]{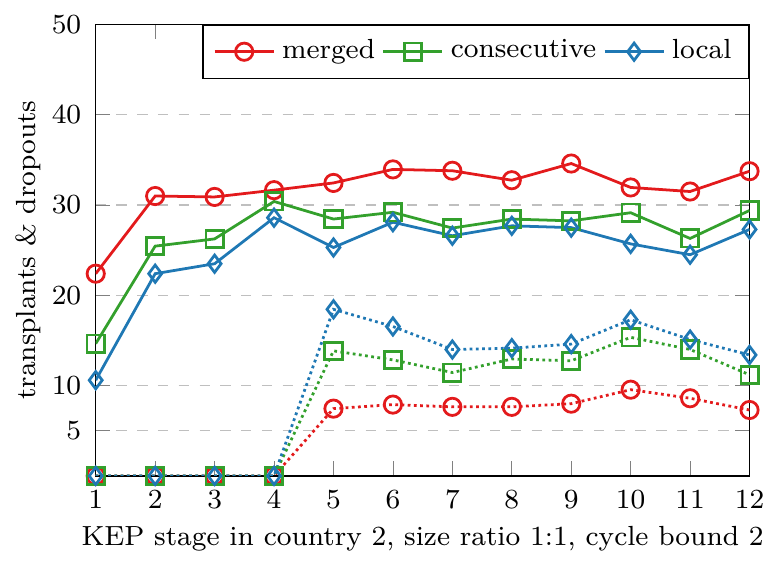}

\includegraphics[width=0.495\textwidth]{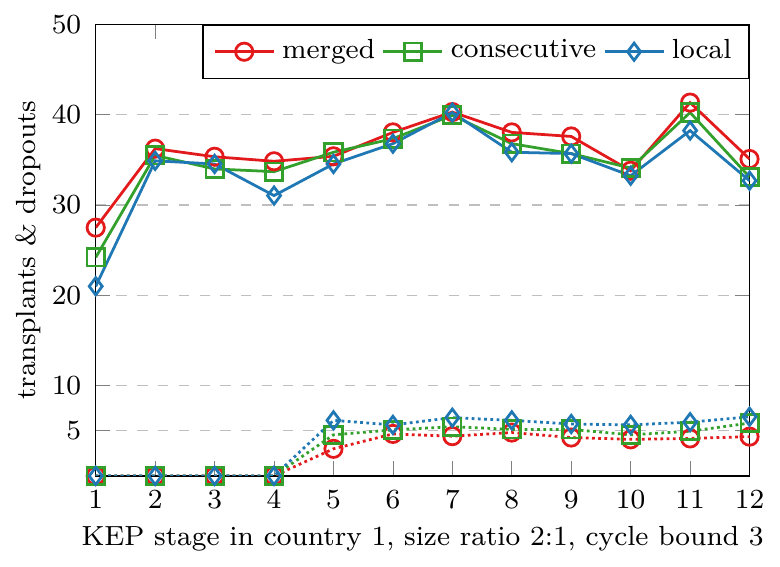}
\includegraphics[width=0.495\textwidth]{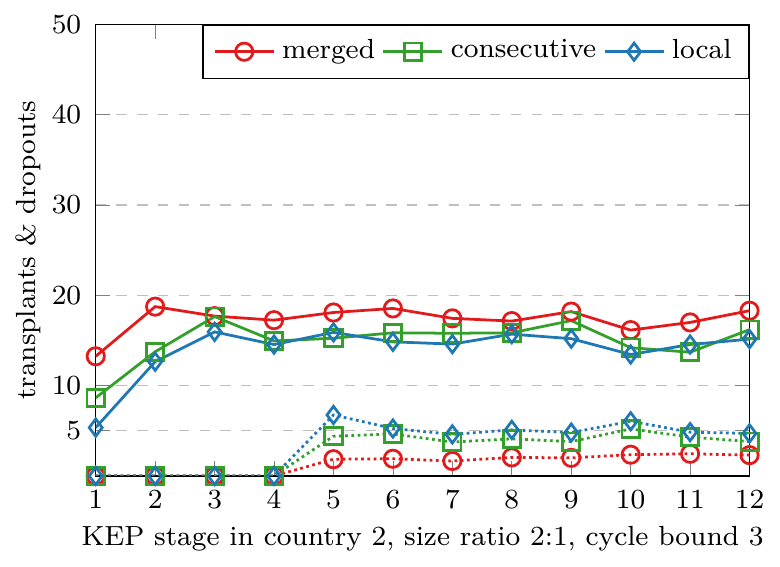}

\includegraphics[width=0.495\textwidth]{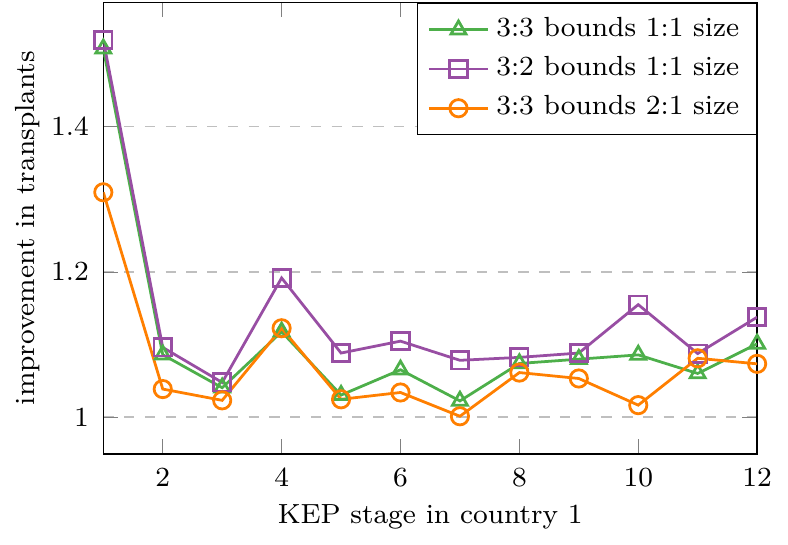}
\includegraphics[width=0.495\textwidth]{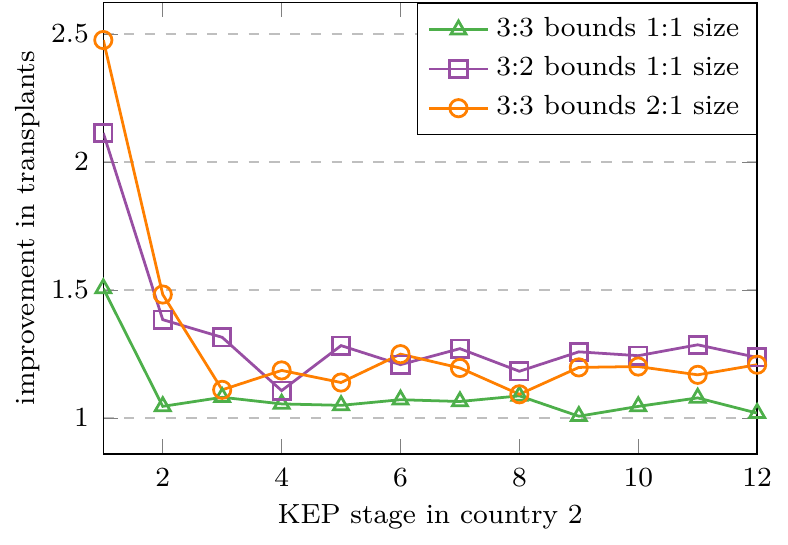}
\end{figure}

Our testing reveals some expected and some unexpected results. The expected results are the following:

\begin{enumerate}
\item On average, when countries collaborate, they do not lose out in terms of total number of transplants after 3 years even without enforcing individual rationality constraints.
\item The merged KEP generally shows much better improvement in the number of transplants than the consecutive KEP  for less restricted countries. The consecutive KEP seems more significant (closer to the merged KEP results) for more restricted cases (see 2:2 bound case in Table \ref{table1:bound}).
\item The size of the KEP pool is a significant factor and positively impacts all forms of collaboration (see Figures \ref{fig:surf2:bound2}-\ref{fig:surf4:bound4} and row 3 in Figure \ref{fig:perstageplots} ).
\item When a smaller country (in terms of KEP pool size) collaborates with a larger country, the smaller country sees greater benefit than the larger country, while the larger country does not lose transplants over what they can achieve by themselves.
\end{enumerate}

\begin{figure}[H]
\caption{C1 and C2 have same pool sizes}
\label{fig:surf0:samepool}
\includegraphics[width=0.495\textwidth]{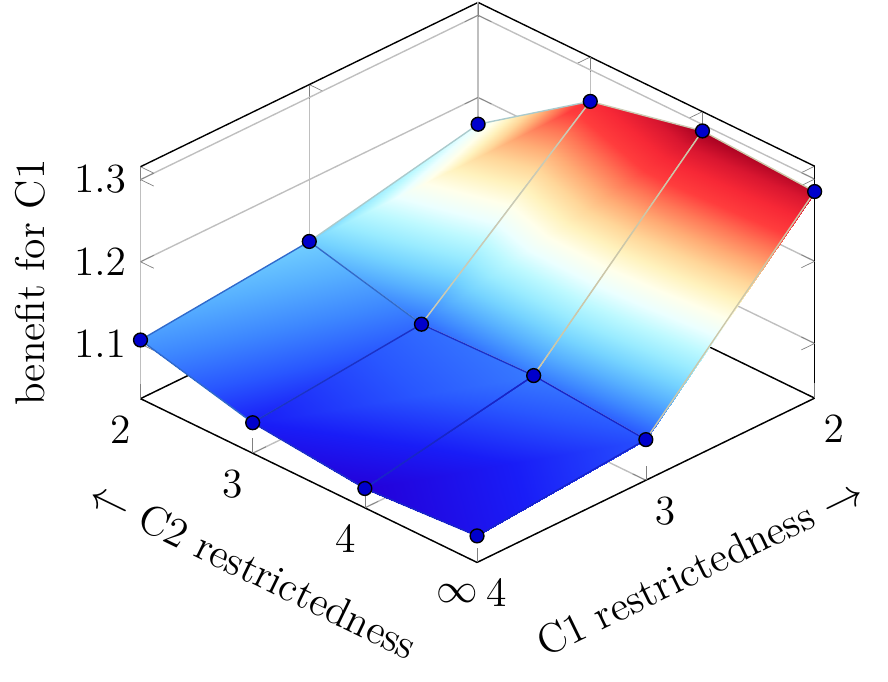}
\includegraphics[width=0.495\textwidth]{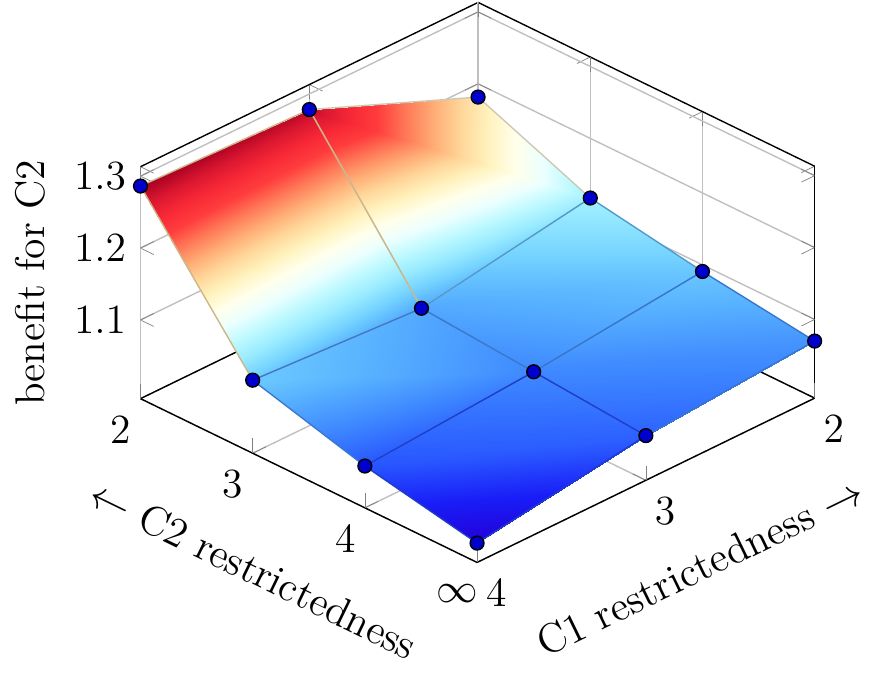}
\end{figure}

\begin{figure}[H]
\caption{C2's pool size is half of C1's pool size ($\frac{C1}{C2}=\frac{2}{1}$)}
\label{fig:surf1:halfpool}
\includegraphics[width=0.495\textwidth]{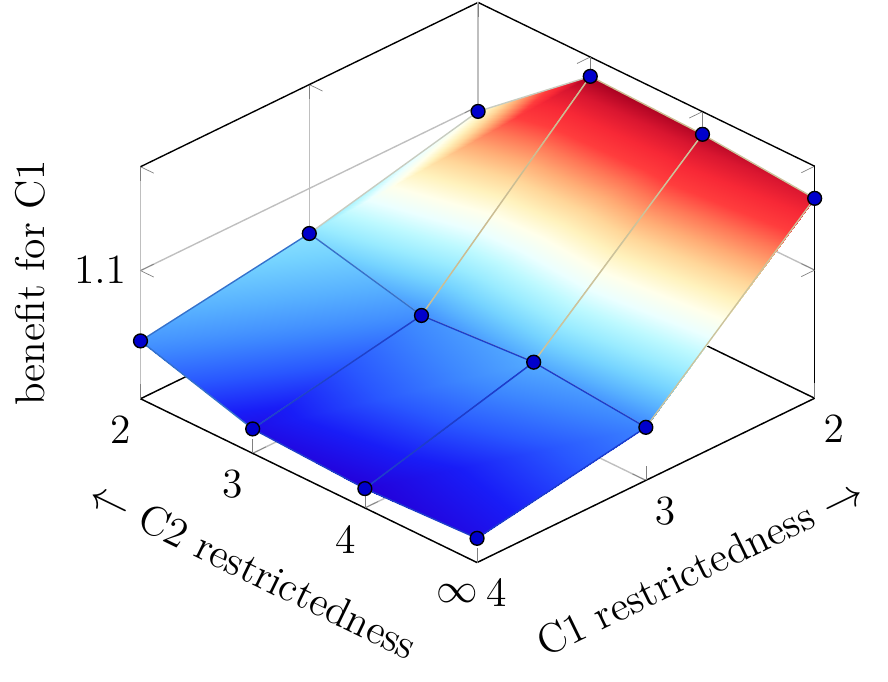}
\includegraphics[width=0.495\textwidth]{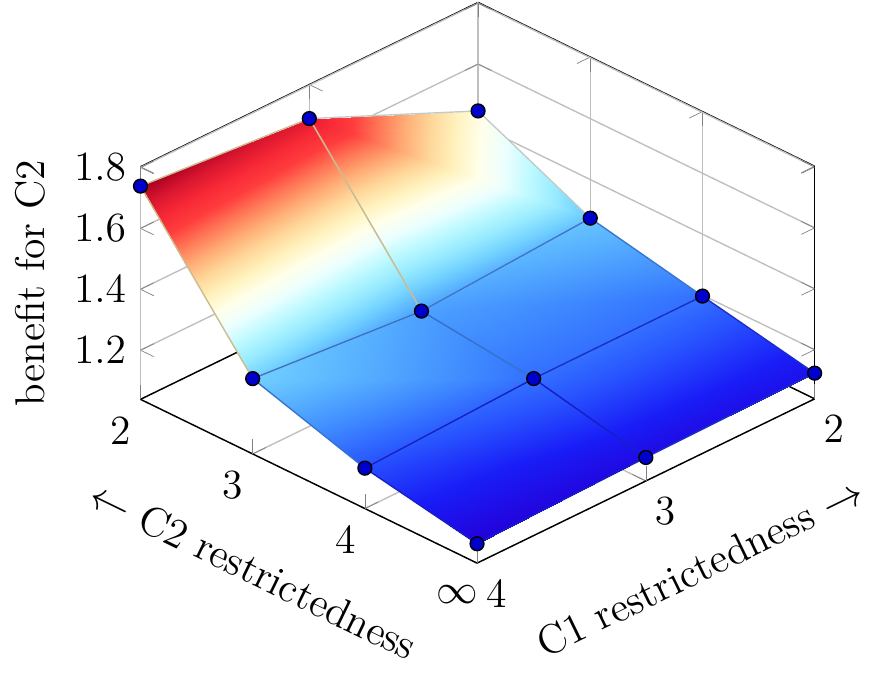}
\end{figure}

Less expected results:
\begin{enumerate}
\item The local restrictions of participating countries have an interesting effect on collaboration (see Figures \ref{fig:surf0:samepool}, \ref{fig:surf1:halfpool}).
\begin{enumerate}
\item[a.] Countries with cycle bound 2 prefer less restricted partners.
\item[b.] Countries with other bounds prefer more restricted partners.
\item[c.] Partners with the same bounds would prefer that their partner had a different bound than themselves. This unusual observation appears consistently as a depression at surface coordinates (2,2), (3,3) and (4,4) in Figures \ref{fig:surf0:samepool} and \ref{fig:surf1:halfpool} and as a valley at bounds 2, 3, 4 in Figures \ref{fig:surf2:bound2}, \ref{fig:surf3:bound3} and \ref{fig:surf4:bound4}, respectively.

\end{enumerate}
\item All things being equal, the partner with the larger KEP pool is preferred. However, there is an exception to this ``larger is better'' rule: sometimes a smaller partner is preferred if they have more suitable restrictions than an equally sized partner as can be seen in Figure \ref{fig:surf3:bound3}: compare the height at $(bound,size)$ coordinates $(2,\frac{3}{4})$ and $(3,\frac{4}{4})$.
\end{enumerate}

We give the following explanation for the anomalous results (especially case 1c): the strange behavior is a result of our chosen policy regarding international transplants. Consider the following example. When countries (C1, C2) have bounds (3, 2), respectively, the increase in the number of transplants for country 1 relative to the (3, 3) bounds case is due to the under-utilisation of the second country's patient-donor pool, leading to a surplus that helps country~1. On the other hand, when (C1, C2) have bounds (3, 4), the increase in the number of transplants relative to the (3, 3) bounds case is caused by the policy that international cycles are of length up to 4 (with at most 2 PDPs in country 1). This policy helps country 1 to secure more transplants for itself since the international cycles are less restricted than the (3, 3) bounds case.


\newpage

\begin{table}[H]
\caption{Number of transplants in each country (average of 20 instances). Same sized pools, different policy constraints.}
\label{table1:bound}

\centering
\begin{tabular}{| c | c | c c c | c c c |}
\hline
C1 & C2 & C1 & C1 & C1 & C2 & C2 & C2\\
bound & bound & local & seq. & merged & local & seq. & merged\\
\hline
2 & 2 & 301.90 & 343.65 & 352.45 & 297.80 & 338.65 & 351.85\\
2 & 3 & 301.90 & 325.15 & 380.90 & 408.40 & 429.95 & 456.20\\
2 & 4 & 301.90 & 315.05 & 390.10 & 433.90 & 445.10 & 473.30\\
2 & $\infty$ & 301.90 & 307.35 & 387.90 & 449.60 & 458.10 & 481.05\\
\hline
3 & 2 & 408.85 & 434.15 & 459.85 & 297.80 & 323.40 & 380.70\\
3 & 3 & 408.85 & 418.80 & 445.85 & 408.40 & 414.60 & 440.25\\
3 & 4 & 408.85 & 414.15 & 447.45 & 433.90 & 436.25 & 462.35\\
3 & $\infty$ & 408.85 & 410.15 & 442.75 & 449.60 & 450.80 & 473.30\\
\hline
4 & 2 & 433.75 & 447.85 & 479.20 & 297.80 & 310.65 & 383.05\\
4 & 3 & 433.75 & 438.25 & 464.35 & 408.40 & 411.65 & 446.05\\
4 & 4 & 433.75 & 436.00 & 458.45 & 433.90 & 433.85 & 455.05\\
4 & $\infty$ & 433.75 & 434.35 & 462.25 & 449.60 & 449.60 & 457.45\\
\hline
\end{tabular}
\end{table}

\begin{table}[H]
\caption{Number of transplants in each country (average of 20 instances). Country 1 has fixed pool size of 100\% and cycle bound 3, Country 2 has varying pool size and varying policy constraints.}
\label{table2:sizebound}

\centering
\begin{tabular}{| c | c | c c c | c c c |}
\hline
C2 & C2 & C1 & C1 & C1 & C2 & C2 & C2\\
bound & size & local & seq. & merged & local & seq. & merged\\
\hline
2 & 1/4 & 408.85 & 425.35 & 425.25 & 39.6 & 60.35 & 95\\
2 & 2/4 & 408.85 & 431.85 & 439.7 & 110.6 & 135.55 & 187\\
2 & 3/4 & 408.85 & 432.6 & 449.4 & 202.5 & 228.1 & 284.55\\
2 & 4/4 & 408.85 & 434.15 & 459.85 & 297.8 & 323.4 & 380.7\\
\hline
3 & 1/4 & 408.85 & 422.1 & 422.5 & 63.3 & 77.7 & 102.1\\
3 & 2/4 & 408.85 & 420.5 & 433.6 & 167.9 & 179.05 & 207.85\\
3 & 3/4 & 408.85 & 419.6 & 440.8 & 289.2 & 298.55 & 325.55\\
3 & 4/4 & 408.85 & 418.8 & 445.85 & 408.4 & 414.6 & 440.25\\
\hline
4 & 1/4 & 408.85 & 421.3 & 423.45 & 72.1 & 85.55 & 109.25\\
4 & 2/4 & 408.85 & 417.7 & 435.35 & 186.65 & 195.6 & 223.15\\
4 & 3/4 & 408.85 & 415.6 & 442.2 & 312 & 316 & 345.1\\
4 & 4/4 & 408.85 & 414.15 & 447.45 & 433.9 & 436.25 & 462.35\\
\hline
$\infty$ & 1/4 & 408.85 & 419.6 & 422.5 & 79.25 & 86.75 & 108.8\\
$\infty$ & 2/4 & 408.85 & 414.9 & 433 & 203.35 & 205.2 & 226.9\\
$\infty$ & 3/4 & 408.85 & 411.7 & 440.45 & 328.95 & 330.55 & 351.3\\
$\infty$ & 4/4 & 408.85 & 410.15 & 442.75 & 449.6 & 450.8 & 473.3\\
\hline
\end{tabular}
\end{table}

\section{Conclusion}

We studied the multi-country kidney exchange problem, where the participating countries may have different constraints and objectives on their national cycles and chains and the parts of the international cycles and chains they are involved in. We formulated IP-models to describe various reasonable conditions by extending the classical cycle- and edge-formulation models with new segment variables and with constraints linking the three types of variables. These formulations are particularly useful if some participant country has no bounds of the lengths of the cycles or altruistic chains. 

In the simulation part we tested the two-country scenario focusing on the dependence of the countries' benefits on the sizes of their pools, and their restrictions. The simulations are realistic with regard to the current practices in Europe and can provide interesting consideration for the decision makers.

We briefly described the way that altruistic chains can be treated in this framework, but in a more extensive study one could consider the particular challenges and constraints related to them. In a potential follow up paper one could also extend this study by considering also the quality of the transplants, rather than only the number of transplants, which is the primary, but not the only goal of the current kidney exchange programmes.  

\section*{Acknowledgement}

The authors acknowledge the financial support of the ENCKEP (European Network for Collaboration on Kidney Exchange Programmes) COST Action given for the short-term scientific missions for Bir\'o, Gyetvai, Mincu and Popa.

\bibliographystyle{plain}
\bibliography{bibliography}

\newpage
\section*{Appendix}
Hereby we describe the cycle-search algorithm used in our simulations. An important fact is that there very little data dependency, meaning that the code is highly parallelizable. In our implementation we simply parallelize by the outermost \textit{for} loop, distributing each $i$ to a different worker thread and then collect the results back to the main thread after each worker thread is finished. We use Java and Java Threads for a simple and portable solution.

\begin{algorithm}
\floatname{CycleSearch}
\caption{\textbf{Algorithm CycleSearch}(\textit{Set V}) adds cycle variables to the model that respect all the constraints. The functions in bold are defined separately. Variable \textit{InternationalCycleLength} is assumed to be a globally available parameter.}
\label{alg:cyclesearch}
\hrule
\begin{algorithmic}[1]
\For {$i \gets 1,|V|$ }
	\For {$j \gets i+1,|V|$}
		\If {$(i,j) \in A \textbf{ and } (j,i) \in A $}
			\If{\textbf{IsValidCycle}$((i,j))$}										
				\State add cycle $(i,j)$ to model
			\EndIf
		\If{$\text{InternationalCycleLength}<3$}
		\State \textbf{continue}
		\EndIf
		\EndIf
		\For {$k \gets j+1,|V|$}
			\If {$(i,j) \in A \textbf{ and } (j,k) \in A \textbf{ and } (k,i) \in A$}
				\If{\textbf{IsValidCycle}$((i,j,k))$}			
					\State add cycle $(i,j,k)$ to model
				\EndIf			
			\EndIf
			\If {$(i,k) \in A \textbf{ and } (k,j) \in A \textbf{ and } (j,i) \in A$}
				\If{\textbf{IsValidCycle}$((i,k,j))$}			
					\State add cycle $(i,k,j)$ to model
				\EndIf			
			\EndIf
					\If{$\text{InternationalCycleLength}<4$}
		\State \textbf{continue}
		\EndIf		
			\For {$\ell \gets k+1,|V|$}
				\If {$(i,j) \in A \textbf{ and } (j,k) \in A \textbf{ and } (k,\ell) \in A \textbf{ and } (\ell,i) \in A$}
					\If{\textbf{IsValidCycle}$((i,j,k,\ell))$}			
						\State add cycle $(i,j,k,\ell)$ to model
					\EndIf			
				\EndIf
				\If {$(i,j) \in A \textbf{ and } (j,\ell) \in A \textbf{ and } (\ell,k) \in A \textbf{ and } (k,i) \in A$}
					\If{\textbf{IsValidCycle}$((i,j,\ell,k))$}			
						\State add cycle $(i,j,\ell,k)$ to model
					\EndIf			
				\EndIf
				
				\If {$(i,k) \in A \textbf{ and } (k,j) \in A \textbf{ and } (j,\ell) \in A \textbf{ and } (\ell,i) \in A$}
					\If{\textbf{IsValidCycle}$((i,k,j,\ell))$}			
						\State add cycle $(i,k,j,\ell)$ to model
					\EndIf			
				\EndIf
				\If {$(i,k) \in A \textbf{ and } (k,\ell) \in A \textbf{ and } (\ell,j) \in A \textbf{ and } (j,i) \in A$}
					\If{\textbf{IsValidCycle}$((i,k,\ell,j))$}			
						\State add cycle $(i,k,\ell,j)$ to model
					\EndIf			
				\EndIf
				
				\If {$(i,\ell) \in A \textbf{ and } (\ell,j) \in A \textbf{ and } (j,k) \in A \textbf{ and } (k,i) \in A$}
					\If{\textbf{IsValidCycle}$((i,\ell,j,k))$}			
						\State add cycle $(i,\ell,j,k)$ to model
					\EndIf			
				\EndIf
				\If {$(i,\ell) \in A \textbf{ and } (\ell,k) \in A \textbf{ and } (k,j) \in A \textbf{ and } (j,i) \in A$}
					\If{\textbf{IsValidCycle}$((i,\ell,k,j))$}			
						\State add cycle $(i,\ell,k,j)$ to model
					\EndIf			
				\EndIf
				
			\EndFor
		\EndFor
	\EndFor
\EndFor
\end{algorithmic}
\end{algorithm}

\begin{algorithm}
\floatname{IsValidCycle}
\caption{\textbf{Function IsValidCycle}(\textit{Cycle} $c$) returns Boolean: whether cycle variable must be added to model.  $\text{MaxCycleLength}_{\text{country}}$ is a globally available parameter. The \textit{Country} function returns the country affiliation of the PDP. Other functions are defined separately.}
\hrule
\begin{algorithmic}[1]
\If{\textbf{IsLocal}$(c)$ \textbf{and} $|c|\leq \text{MaxCycleLength}_{\text{Country}(c_1)}$}
\State \textbf{return} true
\ElsIf {\textbf{CheckCountries}$(c)$ \textbf{and CheckSeg}$(c)$}
\State \textbf{return} true
\EndIf
\State \textbf{return} false
\end{algorithmic}
\end{algorithm}

\begin{algorithm}
\floatname{IsLocal}
\caption{\textbf{Function IsLocal}(\textit{Cycle c}) returns Boolean: whether \textit{Cycle c} is not an international cycle. The \textit{Country} function returns the country affiliation of the PDP.}
\label{alg:islocal}
\hrule
\begin{algorithmic}[1]
\For {$i \gets 2, |c| $ }
	\If{$\text{Country}(c_1) \neq \text{Country}(c_i)$}
	\State \textbf{return} false
	\EndIf
\EndFor
	\State \textbf{return} true
\end{algorithmic}
\end{algorithm}

\begin{algorithm}
\floatname{CheckCountries}
\caption{\textbf{Function CheckCountries}(\textit{Cycle c}) returns Boolean: whether the country constraints are verified for \textit{Cycle c}. Variables that are not initialized here are assumed to be globally available parameters. The \textit{Country} function returns the country affiliation of the PDP.}
\label{alg:checkcountries}
\hrule
\begin{algorithmic}[1]
\For {$\text{country} \gets 1,\text{NumCountries}$ }
	\State $\text{Count}_{\text{country}} \gets 0$
\EndFor
\For {$i \gets 1, |c| $ }
	\State $\text{country} \gets \text{Country}(c_i)$
	\State $\text{Count}_{\text{country}} \gets \text{Count}_{\text{country}} +1$
	\If{$\text{Count}_{\text{country}} > \text{MaxCount}_{\text{country}}$}
	\State \textbf{return} false
	\EndIf
\EndFor
\State $\text{participants} \gets 0$
\For {$\text{country} \gets 1,\text{NumCountries}$ }
	\If{$\text{Count}_{\text{country}} > 0$}
	\State $\text{participants} \gets \text{participants}+1$
	\EndIf
	\If{$\text{participants} > \text{MaxParticipantCountries}$}
	\State \textbf{return} false
	\EndIf
\EndFor
\State \textbf{return} true
\end{algorithmic}
\end{algorithm}

\begin{algorithm}
\floatname{CheckSeg}
\caption{\textbf{Function CheckSeg}(\textit{Cycle c}) returns Boolean: whether the segment constraints are verified for \textit{Cycle c}. Variables that are not initialized here are assumed to be globally available parameters. The \textit{Country} function returns the country affiliation of the PDP.}
\label{alg:checksegments}
\hrule
\begin{algorithmic}[1]
\For {$i \gets 1, |c| $ }
	\State $\text{NodeList.append}(c_i)$
\EndFor
\State $\text{NodeList.append}(c_1)$
\State $\text{seglength} \gets 0$
\For {$\text{country} \gets 1,\text{NumCountries}$ }
	\State $\text{NumSegments}_{\text{country}} \gets 0$
\EndFor
\For {$i \gets 1,\text{length}(\text{NodeList})-1$ }
	\State $\text{source} \gets \text{NodeList}_i$
	\State $\text{target} \gets \text{NodeList}_{i+1}$
	\If {Country(source) = Country(target)}
		\State $\text{seglength} \gets \text{seglength}+1$
		\If {$\text{seglength } > \text{ MaxSegmentLength}_{\text{Country(target)}}$}
		\State \textbf{return} false
		\EndIf
\Else
		\State $\text{seglength} \gets 0$
		\State $\text{NumSegments}_{\text{Country(target)}} \gets \text{NumSegments}_{\text{Country(target)}} +1$
		\If {$\text{NumSegments}_{\text{Country(target)}} > \text{MaxSegments}_{\text{Country(target)}}$}
		\State \textbf{return} false
	\EndIf
	\EndIf		
\EndFor
\end{algorithmic}
\end{algorithm}

\end{document}